\documentclass[11 pt]{article}

\usepackage{amssymb}

\def\nmark{\mbox{$\rm\bf\kern0.2em\rule{0.06em}{1.45ex}\kern-0.3em
N$}}
\def\dmark{\mbox{$\rm\bf\kern0.2em\rule{0.06em}{1.45ex}\kern-0.3em
D$}}
\def\cmark{\mbox{$\rm\bf\kern0.2em\rule{0.06em}{1.45ex}\kern-0.3em
C$}}
\def\rmark{\mbox{$\rm\bf\kern0.2em\rule{0.06em}{1.45ex}\kern-0.3em
R$}}

\begin{document}
\title{\large \bf
Which  weighted composition operators are hyponormal on the Hardy and weighted Bergman spaces?}
 \author{Mahsa Fatehi and Mahmood Haji Shaabani}

{\maketitle}
\begin{abstract}
In this paper, we study hyponormal weighed composition operators on the Hardy and weighted Bergman spaces. For functions $\psi \in A(\mathbb{D})$ which are not the zero function, we  characterize all hyponormal compact weighted composition operators $C_{\psi,\varphi}$ on $H^{2}$ and $A^{2}_{\alpha}$. Next, we show that for $\varphi \in \mbox{LFT}(\mathbb{D})$, if $C_{\varphi}$ is hyponormal on $H^{2}$ or $A^{2}_{\alpha}$, then $\varphi(z)=\lambda z$, where $|\lambda| \leq 1$ or $\varphi$ is a hyperbolic non-automorphism with $\varphi(0)=0$ and such that $\varphi$ has another fixed point in $\partial \mathbb{D}$. After that, we find the essential spectral radius of $C_{\varphi}$ on $H^{2}$ and $A^{2}_{\alpha}$, when $\varphi$ has a Denjoy-Wolff point $\zeta \in \partial \mathbb{D}$. Finally, descriptions of spectral radii are provided for some hyponormal weighted composition operators on $H^{2}$ and $A^{2}_{\alpha}$.
\end{abstract}

\footnote{AMS Subject Classifications. Primary 47B33; Secondary
47B20\\
{\it key words and phrases}:  Hardy space, Weighted Bergman spaces, weighted
composition operator, hyponormal.}

 \section{Introduction}
Let $\mathbb{D}$ denote the open unit disk in the complex plane. The algebra $A(\mathbb{D})$
consists of all continuous functions on the closure of
$\mathbb{D}$ that are analytic on ${\mathbb{D}}$. The Hilbert spaces of primary interest to us will be the Hardy space $H^{2}$ and the weighted Bergman spaces $A_{\alpha}^{2}$. For $f$ which is analytic on $\mathbb{D}$, we denote by $\hat{f}(n)$ the $n$-th coefficient of $f$ in its Maclaurin series. The Hardy space $H^{2}$ is the collection of all such functions $f$ for which
$$\|f\|_{1}^{2}=\sum_{n=0}^{\infty}|\hat{f}(n)|^{2}<\infty.$$
For $\alpha > -1$, the weighted Bergman space $A_{\alpha}^{2}$ consists of all analytic $f$ on $\mathbb{D}$ such that
$$\|f\|_{\alpha+2}^{2}=\int_{\mathbb{D}}|f(z)|^{2}(\alpha+1)(1-|z|^{2})^{\alpha}dA(z)<\infty,$$
where $dA$ is normalized area measure on $\mathbb{D}$. Throughout this paper, let $\gamma=1$ for $H^{2}$ and $\gamma=\alpha+2$ for $A_{\alpha}^{2}$. We know that both the weighted Bergman space and the Hardy space are reproducing kernel Hilbert spaces, when the reproducing kernel for evaluation at $w$ is given by $K_{w}(z)=(1-\overline{w}z)^{-\gamma}$ for $z,w \in \mathbb{D}$. Also the norm of $K_{w}$ is $(1-|w|^{2})^{-\gamma/2}$. We write $H^{\infty}$ for the space of bounded analytic functions on $\mathbb{D}$, with supremum norm  $\|f\|_{\infty}$.\par
Let $\varphi$ be an analytic self-map of $\mathbb{D}$. If $H$ is a Hilbert space of analytic functions on $\mathbb{D}$, the composition operator $C_{\varphi}$ on $H$
 is defined by the rule $C_{\varphi}(f)=f \circ \varphi$. Moreover, for an analytic function $\psi$ on $\mathbb{D}$ and an analytic self-map $\varphi$ of $\mathbb{D}$, we define the weighted composition operator $C_{\psi,\varphi}$ on $H$ by $C_{\psi,\varphi}f=\psi f \circ \varphi$ for all $f \in H$. Such weighted composition operators are clearly bounded on $H^{2}$ and $A_{\alpha}^{2}$ when $\psi$ is bounded on $\mathbb{D}$, but the boundedness of $\psi$ on $\mathbb{D}$ is not necessary for $C_{\psi,\varphi}$ to be bounded. Moreover, if $C_{\psi,\varphi}$ is bounded on $H^{2}$ or $A_{\alpha}^{2}$, then by  \cite[Lemma 1]{Bourdon} and \cite[p. 1211]{fash},
 \begin{eqnarray}
 C_{\psi,\varphi}^{\ast}(K_{w})=\overline{\psi(w)}K_{\varphi(w)}.
  \end{eqnarray}\par
Let $P$ denote the
projection of $L^{2}(\partial \mathbb{D})$ onto $H^{2}$. For each $b \in {L^{\infty}(\partial \mathbb{D})}$,
the Toeplitz  operator $T_{b}$ acts on $H^{2}$ by
$T_{b}(f)=P(bf)$.
Also suppose that $P_{\alpha}$ is the projection of $L^{2}({\mathbb{D}}, dA_{\alpha})$
onto $A^{2}_{\alpha}$.
For each function $w \in L^{\infty}(\mathbb{D})$, the Toeplitz  operator
$T_{w}$ on $A^{2}_{\alpha}$ is defined by $T_{w}(f)=P_{\alpha}(wf)$. Since $P$ and $P_{\alpha}$ are bounded, the Toeplitz  operators are
bounded on $H^{2}$ and $A^{2}_{\alpha}$.\par

A linear-fractional self-map of $\mathbb{D}$ is a map of the form
\begin{eqnarray}
\varphi(z)=\frac{az+b}{cz+d},
\end{eqnarray}
with $ad-bc \neq 0$, for which $\varphi(\mathbb{D}) \subseteq \mathbb{D}$. We denote the set of those maps by $\mbox{LFT}(\mathbb{D})$. It is well known that the automorphisms of $\mathbb{D}$, denoted $\mbox{Aut}(\mathbb{D})$, are the maps in $\mbox{LFT}(\mathbb{D})$ that take $\mathbb{D}$
 onto itself, and that they are of the form $\varphi(z)=\lambda(a-z)/(1-\overline{a}z)$, where $|\lambda|=1$ and $|a| < 1$ (see \cite{c3}). The automorphisms will be classified into three disjoint types:\par
 $\bullet$ Elliptic if one fixed point is in the disk and the other is in the complement of the closed disk.\par
 $\bullet$ Hyperbolic if both fixed points are on the unit circle.\par
 $\bullet$ Parabolic if there is one fixed point on the unit circle of multiplicity $2$.\par

Suppose $\varphi\in\mbox{LFT}(\mathbb{D})$ is
as in Equation (2).
It is well-known that the adjoint of $C_{\varphi}$  acting on $H^{2}$ and $A^{2}_{\alpha}$ is given by
$C_{\varphi}^{\ast}=T_{g}C_{\sigma}T_{h}^{\ast},$
where
$\sigma(z):=({\overline{a}z-\overline{c}})/({-\overline{b}z+\overline{d}})$
is a self-map of $\mathbb{D}$, $g(z):=(-\overline{b}z+\overline{d})^{-\gamma}$ and
$h(z):=(cz+d)^{\gamma}$. Note that  $g$ and $h$ are in  $H^{\infty}$ (see \cite{co2} and \cite{hu}). If $\varphi(\zeta)=\eta$
for $\zeta,\eta \in \partial \mathbb{D}$, then
$\sigma(\eta)=\zeta$. We know that $\varphi$ is an
automorphism if and only if $\sigma$ is, and in this case
$\sigma=\varphi^{-1}$. The map
$\sigma$ is called the Krein adjoint of $\varphi$. We will refer to $g$ and $h$ as the Cowen auxiliary functions for $\varphi$.
From now on, unless otherwise stated, we
assume that $\sigma$, $h$ and $g$ are given as above.
In \cite{kmm} and \cite{mw2}, the adjoint of $C_{\varphi}$, modulo the ideal of compact operators on $H^{2}$ and $A^{2}_{\alpha}$, was obtained, when $\|\varphi\|_{\infty}=1$ but $\varphi$ is not an automorphism; then in \cite{fash} the present  authors gave another proof for the form of this adjoint. \par

We say $\varphi$ has a finite angular derivative at a point $\zeta \in \partial \mathbb{D}$ if there is a point $w \in \partial \mathbb{D}$ such that $(\varphi(z)-w)/(z-\zeta)$ has a finite limit as $z$ tends nontangentially to $\zeta$. This limit, if it exists, is called the angular derivative of $\varphi$ at $\zeta$, and is denoted by $\varphi'(\zeta)$. Throughout this paper, let $F(\varphi)$ denote the set of all points in $\partial \mathbb{D}$ at which $\varphi$ has a finite angular derivative.
\par
Let $\varphi$ be an analytic map from the open unit disk into itself which is neither the identity map nor an elliptic automorphism of $\mathbb{D}$.
The function $\varphi$ can have at most one fixed point inside the open unit disk. If $\varphi$ has a fixed point $a$  inside the open unit disk, then $|\varphi'(a)| \leq 1$. If $\varphi$ has no fixed points inside the open unit disk, then it will have at least one fixed point on the unit circle and for only one of these points, say $a$, $0<\varphi'(a) \leq 1$. The absolute value of the derivative at other fixed points on the unit circle are either greater than $1$ or they do not exist at all. Therefore, it is clear that $\varphi$ has exactly one fixed point $a$ on the closed unit disk satisfying $|\varphi'(a)| \leq 1$. This point is known as the Denjoy-Wolff point of $\varphi$. Also for each $z \in \mathbb{D}$, $\lim_{n \rightarrow \infty}\varphi_{n}(z)=a$, where $\varphi_{0}$ is the identity map and $\varphi_{n}$ denotes the $n$-th iterate of $\varphi$.\par

We say that an operator $A$ on a Hilbert space $H$ is hyponormal if $A^{\ast}A-AA^{\ast} \geq 0$, or equivalently if $\|A^{\ast}f\| \leq \|Af\|$ for all $f \in H$.  Recall that an operator $T$ on a Hilbert space $H$ is said to be
normal if $TT^{\ast}=T^{\ast}T$ on $H$. Also $T$ is unitary if
$TT^{\ast}=T^{\ast}T=I$. The normal composition operators on
$A^{2}_{\alpha}$ and $H^2$ have symbol $\varphi(z)=az$, where $|a|\leq
1$ (see \cite[Theorem 8.2]{cm1}).
 Also, it is easy to see that only
the rotation maps $\varphi(z)=\zeta z$, $|\zeta|=1$, induce
unitary composition operators $C_{\varphi}$ on $H^{2}$ and
$A^{2}_{\alpha}$. The normal and unitary weighted composition
operators on $H^{2}$ and $A^{2}_{\alpha}$ were investigated in \cite{Bourdon} and  \cite{l}. Cowen et al. \cite{ck} found the relationship between properties of the symbol $\varphi$ and the hyponormality of composition operators $C_{\varphi}$.
After that Zorboska \cite{z} investigated the hyponormal composition operators on the weighted Hardy spaces. Recently in \cite{coko}, Cowen et al. investigated that when $C_{\psi,\varphi}^{\ast}$ is hyponormal. In this paper, we are interested in the hyponormal composition operators and weighted composition operators on $H^{2}$ and $A_{\alpha}^{2}$.  Our description of all hyponormal compact weighted composition operators $C_{\psi,\varphi}$ on $H^{2}$ and $A_{\alpha}^{2}$ induced by nonzero functions $\psi \in A(\mathbb{D})$, which appears in Section 2, reveals that for such operators both $\psi$ and $\varphi$ are linear-fractional self-maps. Moreover, in Section 2, we see that for $\varphi \in \mbox{LFT}(\mathbb{D})$, if $C_{\varphi}$ is hyponormal on $H^{2}$ or $A_{\alpha}^{2}$, then $C_{\varphi}$ is normal or $\varphi$ is a hyperbolic non-automorphism with $\varphi(0)=0$ and such that $\varphi$ has another fixed point in $\partial \mathbb{D}$. In the final section of this paper, we investigate the spectral radii of some hyponormal weighted composition operators on $H^{2}$ and $A_{\alpha}^{2}$.

\section{Hyponormal weighted composition operators on $H^{2}$ and $A^{2}_{\alpha}$}

In this section, we find some hyponormal weighted composition operators on $H^{2}$ and $A_{\alpha}^{2}$. Let $H$ be a Hilbert space. The set of all bounded operators and the set of
all compact operators from $H$ into itself are denoted by $B(H)$
and $B_{0}(H)$, respectively.
Suppose that $T$ belongs to $B(H^{2})$ or $B(A_{\alpha}^{2})$. Through this paper, the spectrum of $T$, the essential spectrum of $T$, the approximate point spectrum of $T$ and the point spectrum of $T$ are denoted by $\sigma_{\gamma}(T)$, $\sigma_{e,\gamma}(T)$, $\sigma_{ap,\gamma}(T)$ and $\sigma_{p,\gamma}(T)$, respectively, for $H^{2}$ and $A_{\alpha}^{2}$. Also the spectral radius of $T$ and the essential spectral radius of $T$ are denoted by $r_{\gamma}(T)$ and $r_{e,\gamma}(T)$, respectively. Moreover, we denote by $\|T\|_{\gamma}$ and $\|T\|_{e,\gamma}$ the norm of the operator $T$ and the essential norm of the operator $T$, respectively, on $H^{2}$ and $A_{\alpha}^{2}$. In this paper, for $r \in \mathbb{C}$ and a set $A \subseteq \mathbb{C}$, we define $rA:=\{ ra:a\in A\}$. \par
In \cite{Bourdon1}, Bourdon proved the next lemma for $H^{2}$. We show that \cite[Lemma 5.1]{Bourdon1} holds for $A^{2}_{\alpha}$. The idea of the proof of the following lemma is inspired by \cite[Lemma 5.1]{Bourdon1}. \\ \par

{\bf Lemma 2.1.} {\it Suppose that $\varphi$, not the identity and not an elliptic automorphism of $\mathbb{D}$, is an analytic map of the unit disk into itself with Denjoy-Wolff point $\zeta$. Assume that $\psi \in H^{\infty}$
extends to be continuous on $\mathbb{D} \cup \{\zeta\}$ (if $\zeta \in \partial \mathbb{D}$). Suppose that $C_{\psi,\varphi}$ is considered as an operator on $A_{\alpha}^{2}$. If $\lambda$ is an eigenvalue of $C_{\psi,\varphi}$, then $|\lambda| \leq |\psi(\zeta)|r_{\alpha+2}(C_{\varphi})$. If $\psi(\zeta)=0$
and $\varphi$ and $\psi$ are nonconstant, then $C_{\psi,\varphi}$ has no eigenvalues.}\bigskip

{\bf Proof.} Let $\lambda$ be an eigenvalue for $C_{\psi,\varphi}$ with corresponding eigenvector $f \in A_{\alpha}^{2}$. Since $C_{\psi,\varphi}^{n}f=\lambda^{n}f$, for each positive integer $n$ and $z \in \mathbb{D}$, we have
\begin{eqnarray}
 f(\varphi_{n}(z)) \prod_{j=0}^{n-1}\psi(\varphi_{j}(z))=\lambda^{n}f(z).
  \end{eqnarray}
For any fixed point $z \in \mathbb{D}$ and positive integer $n$, one can easily see that
\begin{eqnarray}
|f(\varphi_{n}(z))|
&=&|\langle f \circ \varphi_{n},K_{z}\rangle |\nonumber\\
 & \leq & \|f \circ \varphi_{n}\|_{\alpha+2}\|K_{z}\|_{\alpha+2}\nonumber\\
& \leq & \|C_{\varphi}^{n}\|_{\alpha+2}\frac{\|f\|_{\alpha+2}}{(1-|z|^{2})^{(\alpha+2)/2}}.
\end{eqnarray}
Because $f$ is not the zero function, there is $z \in \mathbb{D}$  such that $f(z) \neq 0$. Since $\{z\}$ is a compact set, by the Denjoy-Wolff Theorem,
$\varphi_{j}(z) \rightarrow \zeta$ as $j \rightarrow \infty$. Hence $\psi(\varphi_{j}(z)) \rightarrow \psi(\zeta)$ as $j \rightarrow \infty$.
Equations (3) and (4) yield
\begin{eqnarray}
|\lambda^{n}f(z)|^{1/n}
&\leq& \left|\prod_{j=0}^{n-1}\psi (\varphi_{j}(z))\right|^{1/n}\left(\frac{\|f
\|_{\alpha+2}\|{C_{\varphi}^{n}}\|_{\alpha+2}}{(1-|z|^{2})^{(\alpha+2)/2}}\right)^{1/n}.
\end{eqnarray}
Letting $n \rightarrow \infty$,  by Equation (5) we have
\begin{eqnarray}
|\lambda| \leq |\psi(\zeta)|r_{\alpha+2}(C_{\varphi}).
\end{eqnarray}
Now suppose that $\psi(\zeta)=0$. Hence by Equation (6), $\sigma_{p,\alpha+2}(C_{\psi,\varphi}) \subseteq \{0\}$. If $0 \in \sigma_{p,\alpha+2}(C_{\psi,\varphi})$, then there is $f \in A_{\alpha}^{2}$ such that $f$ is an eigenvector for $C_{\psi,\varphi}$ with corresponding eigenvalue $0$ and $\psi f \circ \varphi \equiv 0$. Since $\psi$ is not the zero function and $\varphi$ is not a constant function, the Open Mapping Theorem implies that $f \equiv 0$ and it is a contradiction.\hfill $\Box$ \\ \par

{\bf Proposition 2.2.} {\it Suppose that $\varphi$, not the identity and not an elliptic automorphism of $\mathbb{D}$, is an analytic map of the unit disk into itself with Denjoy-Wolff point $\zeta$. Assume that $\psi \in H^{\infty}$
extends to be continuous on $\mathbb{D} \cup \{\zeta\}$ (if $\zeta \in \partial \mathbb{D}$). Suppose that $\varphi$ and $\psi$ are nonconstant. If $\psi(\zeta)=0$ and $\sigma_{e,\gamma}(C_{\psi,\varphi})=\{0\}$, then $C_{\psi,\varphi}$ is not hyponormal on $H^{2}$ and $A_{\alpha}^{2}$.}
\bigskip

{\bf Proof.} Assume that $C_{\psi,\varphi}$ is hyponormal on $H^{2}$ and $A_{\alpha}^{2}$. Since $\psi(\zeta)=0$, by \cite[Lemma 5.1]{Bourdon1} and Lemma 2.1, $\sigma_{p,\gamma}(C_{\psi,\varphi})=\emptyset$. Invoking \cite[Proposition 6.7, p. 210]{c1} and \cite[Proposition 4.4, p. 359]{c1}, we see that $\partial \sigma_{\gamma}(C_{\psi,\varphi})=\{0\}$. Thus, $r_{\gamma}(C_{\psi,\varphi})=0$. Then \cite[Theorem 1]{stamp1} ensures that $\|C_{\psi,\varphi}\|_{\gamma}=0$. Hence $\psi\equiv0$ and it is a contradiction.\hfill $\Box$ \\ \par

{\bf Corollary 2.3.} {\it Let $\varphi$ be an analytic self-map of $\mathbb{D}$ with Denjoy-Wolff point $\zeta$, where $|\zeta|=1$. Assume that $\psi \in H^{\infty}$ is continuous at $\zeta$. Then $C_{\psi,\varphi}$ is not a hyponormal compact operator on $H^2$ and $A_{\alpha}^{2}$.}\bigskip

{\bf Proof.} Assume that $C_{\psi,\varphi}$ is a hyponormal compact operator on $H^{2}$ or $A_{\alpha}^{2}$. If $\psi$ is a constant function, then $C_{\varphi}$ is hyponormal. Thus, by \cite[Theorem 1]{z}, $\varphi(0)=0$ and it is a contradiction. Now assume that $\psi$ is nonconstant.
 We have $\sigma_{e,\gamma}(C_{\psi,\varphi})=\{0\}$ and $\|C_{\psi,\varphi}\|_{e,\gamma}=0$. If $Q$ is an arbitrary compact operator on $H^{2}$ or $A_{\alpha}^{2}$, then by \cite[Theorem 2.17]{cm1}, Equation (1) and the Julia-Carath\'eodory Theorem, we have
\begin{eqnarray*}
\|C_{\psi,\varphi}-Q\|_{\gamma}
&\geq & \lim_{r \rightarrow 1}\left\|(C_{\psi,\varphi}-Q)^{\ast}\frac{K_{r\zeta}}{\|K_{r\zeta}\|_{\gamma}}\right\|_{\gamma}\\
&=&\lim_{r \rightarrow 1}|\psi(r \zeta)|\frac{\|K_{\varphi(r\zeta)}\|_{\gamma}}{\|K_{r\zeta}\|_{\gamma}}\\
&=&\lim_{r \rightarrow 1}|\psi(r \zeta)|\left(\frac{1-|\varphi(r\zeta)|^{2}}{1-|r \zeta|^{2}}\right)^{-\gamma/2}\\
&=&|\psi( \zeta)||\varphi'(\zeta)|^{-\gamma/2}.\\
  \end{eqnarray*}
 Hence $\|C_{\psi,\varphi}\|_{e,\gamma} \geq |\psi(\zeta)||\varphi'(\zeta)|^{-\gamma/2}$ and so $\psi(\zeta)=0$.
Proposition 2.2 implies that $C_{\psi,\varphi}$ is not hyponormal and it is a contradiction.
\hfill $\Box$ \\ \par


Suppose $\psi$ is not the zero function and $\psi \in A(\mathbb{D})$. In the next theorem, we see that all hyponormal compact weighted composition operators $C_{\psi,\varphi}$ on $H^2$ and $A_{\alpha}^{2}$ were exhibited in \cite[Theorem 10]{Bourdon} and \cite[Theorem 4.3]{l}. \\ \par

{\bf Theorem 2.4.} {\it Let $\varphi$ be an analytic self-map of $\mathbb{D}$. Assume that $\psi$ is not the zero function and $\psi \in A(\mathbb{D})$. The weighted composition operator
$C_{\psi,\varphi}$ is a hyponormal compact operator on $H^2$ or $A_{\alpha}^{2}$ if and only if
 $\psi=\psi(p)\frac{K_{p}}{K_{p}\circ \varphi}$
 and
 $\varphi=\alpha_{p}\circ(\delta\alpha_{p}),$
 where $p \in \mathbb{D}$ is a fixed point of $\varphi$, $\alpha_{p}(z)=(p-z)/(1-\overline{p}z)$ and $\delta$ is a constant such that $|\delta|<1$.}\bigskip

{\bf Proof.} Let $C_{\psi,\varphi}$ be a hyponormal compact operator. Assume that $\psi$ is a constant function. By
\cite[Proposition 4.3, p. 47]{c4}, it is not hard to see that $C_{\varphi}$ is hyponormal. Hence \cite[Corollary 4.9, p. 48]{c4}
implies that $C_{\varphi}$ is normal. Since $C_{\varphi}$ is compact and normal, we see that $\varphi(z)=\delta z$, where $|\delta|<1$. Thus, the result follows. Now assume that $\psi$ is not a constant function. We claim that $\varphi$ has a fixed point in $\mathbb{D}$. Assume that $\varphi$ has no fixed point on $\mathbb{D}$ and so there is $\zeta \in \partial \mathbb{D}$ which is the Denjoy-Wolff point of $\varphi$. By the proof of Corollary 2.3, $\psi(\zeta)=0$. Lemma 2.1 and \cite[Lemma 5.1]{Bourdon1} imply that $C_{\psi,\varphi}$ has no eigenvalues, so by \cite[Theorem 7.1, p. 214]{c1}, $\sigma_{\gamma}(C_{\psi,\varphi})=\{0\}$. Since $\sigma_{e,\gamma}(C_{\psi,\varphi}) \subseteq \sigma_{\gamma}(C_{\psi,\varphi})  $, $\sigma_{e,\gamma}(C_{\psi,\varphi})=\{0\}$. Therefore, by Proposition 2.2, $C_{\psi,\varphi}$ is not hyponormal and it is a contradiction. Hence we conclude that $\varphi$ must have a fixed point $p \in \mathbb{D}$. Since $C_{\psi,\varphi}$ is a hyponormal compact operator, by
\cite[Corollary 4.9, p. 48]{c4}, $C_{\psi,\varphi}$ is normal. Then by \cite[Theorem 10]{Bourdon} and \cite[Theorem 4.3]{l}, $\psi=\psi(p)\frac{K_{p}}{K_{p}\circ \varphi}$ and $\varphi=\alpha_{p}\circ(\delta\alpha_{p})$, where $|\delta| \leq 1$. Assume that $|\delta| = 1$. It is not hard to see that $\alpha_{p}\circ(\delta\alpha_{p})$
is an automorphism of $\mathbb{D}$, when $|\delta|=1$.  Also \cite[Proposition 3.9]{fash} and \cite[Theorem 10]{ma} imply that $C_{\psi,\varphi}$ is not compact and it is a contradiction. Therefore, $|\delta| < 1$.\par
Conversely, it is easy to see that $\overline{\varphi(\mathbb{D})} \subseteq \mathbb{D}$. Hence by \cite[p. 129]{cm1}, $C_{\psi,\varphi}$ is compact. The conclusion follows from \cite[Theorem 10]{Bourdon} and \cite[Theorem 4.3]{l}.
\hfill $\Box$ \\ \par

By the proof of Theorem 2.4, Corollary 2.5 follows immediately. \\ \par

{\bf Corollary 2.5.} {\it Let $C_{\varphi}$ be a hyponormal compact operator on $H^2$ or $A^{2}_{\alpha}$. Then $\varphi(z)=cz$, where $|c|<1$.}\bigskip

{\bf Proposition 2.6.} {\it Let $\varphi$ be an analytic self-map of $\mathbb{D}$.
Assume that $C_{\varphi}$ is hyponormal on $H^2$ or $A^{2}_{\alpha}$. If $\varphi$
has a finite angular derivative at $\zeta$, then $|\varphi'(\zeta)| \geq1$.
 }\bigskip

{\bf Proof.} Assume that $\varphi$ has a finite angular derivative at $\zeta$. Since $C_{\varphi}$ is hyponormal, by \cite[Theorem 1]{z} and Schwarz's Lemma, $(1-|\varphi(w)|)/(1-|w|)\geq1$. Hence
$\liminf_{|w|\rightarrow 1}\frac{1-|\varphi(w)|}{1-|w|}\geq1.$
Therefore, by \cite[Proposition 2.46]{cm1}, $|\varphi'(\zeta)| \geq1$.\hfill $\Box$ \\ \par

{\bf Proposition 2.7.} {\it Suppose that $\varphi$ is  an analytic self-map of $\mathbb{D}$. Assume that $\varphi \in A(\mathbb{D})$ and $\{e^{i\theta}:|\varphi(e^{i\theta})|=1\}$ has only one element $\zeta$.  Let $\psi \in H^{\infty}$ be continuous at $\zeta$. If $\psi(\zeta)=0$, then $C_{\psi,\varphi}$ is not hyponormal on $H^{2}$ and $A_{\alpha}^{2}$. }\bigskip

{\bf Proof.}  Assume that $C_{\psi,\varphi}$ is hyponormal. Let $\psi(\zeta)=0$. Then by \cite[Proposition 2.3]{fash} and \cite[Corollary 2.2]{kmm}, $C_{\psi,\varphi}$ is compact.  We claim that $\varphi$ has a fixed point in $\mathbb{D}$.  Assume that $\varphi$ has no fixed point in $\mathbb{D}$, so  $\zeta$ is the Denjoy-Wolff point of $\varphi$.  Corollary 2.3 shows that $C_{\psi,\varphi}$ is not a hyponormal compact operator and it is a contradiction. Then $\varphi$ must have a fixed point in $\mathbb{D}$. Since $C_{\psi,\varphi}$ is normal, by \cite[Theorem 10]{Bourdon} and \cite[Theorem 4.3]{l}, we can see that $\psi(\zeta) \neq 0$. It is a contradiction. \hfill $\Box$ \\ \par

Suppose that $\varphi \in \mbox{LFT}(\mathbb{D})$ and $\overline{\varphi(\mathbb{D})} \subseteq \mathbb{D}$. Therefore, $\varphi$ has a fixed point in $\mathbb{D}$ and $C_{\psi,\varphi}$ is compact. If $C_{\psi,\varphi}$ is a hyponormal operator, then it is normal. These operators were characterized in \cite[Theorem 10]{Bourdon} and \cite[Theorem 4.3]{l}. \par

Suppose that $\varphi \in \mbox{LFT}(\mathbb{D})$ is not an automorphism but has $\|\varphi\|_{\infty}=1$. We classify $\varphi$ as follows:\par
$\bullet$ Hyperbolic non-automorphism of $\mathbb{D}$ which has a fixed point in $\partial \mathbb{D}$ of multiplicity $1$. Also it has another fixed point in the complement of $\partial \mathbb{D}$.\par
$\bullet$ Parabolic non-automorphism of $\mathbb{D}$ with a fixed point in $\partial \mathbb{D}$ of multiplicity two.\par
$\bullet$ Non-automorphism with sup-norm equal to $1$ such that it does not have a fixed point in $\partial \mathbb{D}$. It necessarily has a fixed point in $\mathbb{D}$. \par
In \cite[Theorem 23]{derek}, Cowen et al. stated that if $\varphi$ is a hyperbolic non-automorphism with Denjoy-Wolff point $\zeta$ in $\partial \mathbb{D}$, then there is no $\psi \in H^{\infty}$ continuous at $\zeta$  such that $C_{\psi,\varphi}$ is hyponormal on $H^{2}$. Now suppose that $\varphi$ is a linear-fractional self-map of $\mathbb{D}$, not an automorphism, which satisfies $\varphi(\zeta)=\eta$, where $\zeta,\eta \in  \partial \mathbb{D}$ and $\zeta \neq \eta$. Also let $\psi \in H^{\infty}$ be continuous at $\zeta$. In Proposition 2.8, we show that $C_{\psi,\varphi}$ is not hyponormal on $H^{2}$ and $A_{\alpha}^{2}$.\par
 If $A$ and $B$ are subsets of complex plane and $\lambda \neq 0$ is  constant such that $A=\lambda B$, then it is not hard to see that $\sup \{|a|: a \in A\}= \lambda \sup \{|b|: b \in B\}$. We use this fact in the proof of the following proposition.\\ \par

{\bf Proposition 2.8.} {\it Suppose that $\varphi$ is a linear-fractional self-map of $\mathbb{D}$, not an automorphism, which satisfies $\varphi(\zeta)=\eta$, where $\zeta,\eta \in  \partial \mathbb{D}$ and $\zeta \neq \eta$. Let $\psi \in H^{\infty}$ be continuous at $\zeta$. Suppose that $\psi$ is not the zero function. Then $C_{\psi,\varphi}$ is not hyponormal on $H^{2}$ or $A_{\alpha}^{2}$. }\bigskip

{\bf Proof.}
By Proposition 2.7,  assume that $\psi(\zeta) \neq 0$.
We claim that there is $t \in (-\infty,-|\psi(\zeta)|^{2}|\varphi'(\zeta)|^{-\gamma}]$ such that $t \in \sigma_{e,\gamma}(C_{\psi,\varphi}^{\ast}C_{\psi,\varphi}-C_{\psi,\varphi}C_{\psi,\varphi}^{\ast})$. This may be seen as follows.  By   \cite[Proposition 2.3]{fash} and \cite[Corollary 2.2]{kmm}, $\sigma_{e,\gamma}(C_{\psi,\varphi}C_{\psi,\varphi}^{\ast})
=|\psi(\zeta)|^{2}\sigma_{e,\gamma}(C_{\varphi}C_{\varphi}^{\ast})$. Now using \cite[Theorem 3.1]{kmm}, \cite[Proposition 3.6]{kmm} and
\cite[Theorem 3.2]{mw2}, we have $\sigma_{e,\gamma}(C_{\psi,\varphi}C_{\psi,\varphi}^{\ast})=|\psi(\zeta)|^{2}|\varphi'(\zeta)|^{-\gamma}
\sigma_{e,\gamma}(C_{\sigma \circ \varphi})$. By the fact which was stated before Proposition 2.8, we can see that
$r_{e,\gamma}(C_{\psi,\varphi}C_{\psi,\varphi}^{\ast})
=|\psi(\zeta)|^{2}|\varphi'(\zeta)|^{-\gamma}r_{e,\gamma}(C_{\sigma \circ \varphi})$. Therefore,
$r_{e,\gamma}(C_{\psi,\varphi}C_{\psi,\varphi}^{\ast})=|\psi(\zeta)|^{2} |\varphi'(\zeta)|^{-\gamma} \lim_{n \rightarrow \infty} \|C_{\sigma \circ  \varphi}^{n}\|_{e,\gamma}^{1/n}$.
  Also \cite[Proposition 3.4]{kmm} implies that $(\sigma\circ\varphi)'(\zeta)=1$. Hence by \cite[Theorem 5.2]{km}, one can easily see that
  $r_{e,\gamma}(C_{\psi,\varphi}C_{\psi,\varphi}^{\ast}) \geq |\psi(\zeta)|^{2}|\varphi'(\zeta)|^{-\gamma}.$
  Since $C_{\psi,\varphi}C_{\psi,\varphi}^{\ast}$ is positive, there is $\lambda \geq |\psi(\zeta)|^{2}|\varphi'(\zeta)|^{-\gamma}$ such that $\lambda \in \sigma_{e,\gamma}(C_{\psi,\varphi}C_{\psi,\varphi}^{\ast})$. Let $z, z^{\ast}$ and $e$ be the cosets of $C_{\psi,\varphi}$,
$C_{\psi,\varphi}^{\ast}$ and $I$ in the Calkin algebra
$B(H^{2})/B_{0}(H^{2})$ or $B(A_{\alpha}^{2})/{B_{0}(A_{\alpha}^{2})}$.
By \cite[Proposition 2.3]{fash}, \cite[Theorem 3.2]{mw2}, \cite[Corollary 2.2]{kmm} and \cite[Theorem 3.1]{kmm}, it is not hard to see that $(zz^{\ast})(z^{\ast}z)=(z^{\ast}z)(zz^{\ast})=0$ and so the $C^{\ast}$-algebra generated by $zz^{\ast}$, $z^{\ast}z$ and
$e$ is commutative. Let $\mathcal{A}$ be the commutative
$C^{\ast}$-algebra generated by $zz^{\ast}$,$z^{\ast}z$ and
$e$. We denote the collection of all nonzero
homomorphisms of $\mathcal{A} \rightarrow \mathbb{C}$ by $\Sigma$. Since $(zz^{\ast})(z^{\ast}z)=(z^{\ast}z)(zz^{\ast})=0$, $m(z^{\ast}z)=0$ or $m(zz^{\ast})=0$ for each $m \in \Sigma$. Therefore, by \cite[Theorem 8.6, p. 219]{c1}, we have $-\lambda \in \sigma_{e,\gamma}(C_{\psi,\varphi}^{\ast}C_{\psi,\varphi}-C_{\psi,\varphi}C_{\psi,\varphi}^{\ast})$. Hence $C_{\psi,\varphi}$ is not hyponormal on  $H^{2}$
and $A_{\alpha}^{2}$.\hfill $\Box$ \\ \par

A map $\varphi \in \mbox{LFT}(\mathbb{D})$ is called parabolic if it has a fixed point $\zeta \in \partial \mathbb{D}$ of multiplicity $2$. The map $\tau(z):=(1+\overline{\zeta}z)/(1-\overline{\zeta}z)$ takes the unit disk onto the right half-plane $\Pi$ and sends $\zeta$ to $\infty$. Therefore, $\phi:=\tau\circ\varphi\circ\tau^{-1}$ is a self-map of $\Pi$ which fixes only $\infty$, and so must be the mapping of translation by some number $t$, where necessarily $\mbox{Re}t \geq 0$.
Note that  $\mbox{Re}t=0$ if and only if $\varphi \in \mbox{Aut}(\mathbb{D})$.  When the number $t$ is strictly positive, we call $\varphi$ a positive parabolic non-automorphism. Among the linear-fractional self-maps of $\mathbb{D}$ fixing $\zeta \in \partial \mathbb{D}$, the parabolic ones are characterized by $\varphi'(\zeta)=1$. \par

Suppose that $\varphi$ is a linear-fractional self-map of $ \mathbb{D}$. In the following theorem, we see that if $C_{\varphi}$ is hyponormal, then $\varphi(z)=\lambda z$, when $|\lambda| \leq 1$ or $\varphi$ is a hyperbolic non-automorphism with a fixed point $\zeta \in \partial \mathbb{D}$.\\ \par

{\bf Theorem 2.9.} {\it Let $\varphi$ be a linear-fractional self-map of $ \mathbb{D}$. If $C_{\varphi}$ is hyponormal on $H^{2}$ or $A_{\alpha}^{2}$, then $\varphi(z)=\lambda z$ or $\varphi(z)=(1-|c|)z/(cz+1)$, where $|\lambda| \leq 1$ and $0 < |c| < 1$.}\bigskip

{\bf Proof.} Let $\varphi \in \mbox{LFT}(\mathbb{D})$ and $C_{\varphi}$ be hyponormal. By \cite{z}, we know that $\varphi(0)=0$.
We break the proof into five parts.\par
(1) Let $\varphi \in \mbox{Aut}(\mathbb{D})$. Since $\varphi(0)=0$, we can easily see that $\varphi(z)=\lambda z$, where $|\lambda|=1$.\par
(2) Let $\overline{\varphi(\mathbb{D})} \subseteq \mathbb{D}$. Then by \cite[p. 129]{cm1}, $C_{\varphi}$ is compact. \cite[Corollary 4.9, p. 48]{c4} implies that $C_{\varphi}$ is normal. Hence $\varphi(z)=\lambda z$, where $|\lambda| < 1$.\par
(3) Suppose that $\varphi$ is a linear-fractional self-map of $\mathbb{D}$, not an automorphism, which satisfies $\varphi(\zeta)=\eta$ for some  $\zeta,\eta \in \partial \mathbb{D}$ and $\zeta \neq \eta$. Then by Proposition 2.8, $C_{\varphi}$ is not hyponormal.\par
(4) Let $\varphi$ be a parabolic non-automorphism of $\mathbb{D}$. Then $\varphi$ has only one fixed point in $\partial \mathbb{D}$ and it is a contradiction.\par
(5) Let $\varphi$ be a hyperbolic non-automorphism with a fixed point $\zeta \in \partial \mathbb{D}$. Since $\varphi(0)=0$, we can easily see that $\varphi$ is a map of the form $\varphi(z)=az/(cz+1)$, where $a,c \in \mathbb{C}$ and $c \neq 0$. We have $\zeta=a \zeta/(c\zeta+1)$, so $a=c\zeta+1$. If $\sigma$ is the Krein adjoint of $\varphi$, then $\sigma(z)=(\overline{c\zeta}+1)z-\overline{c}$. We have $|c| < 1$, because $\sigma(0) \in \mathbb{D}$.  By \cite[Proposition 3.4]{kmm}, $\varphi \circ \sigma$ is a positive parabolic non-automorphism. Let $\tau(z)=(1+\overline{\zeta }z)/(1-\overline{\zeta}z)$. Then $\tau  \circ \varphi \circ \sigma \circ \tau^{-1}$ sends $0$ to a strictly positive number. One sees, after some patient calculation, that $\tau(\varphi(\sigma(\tau^{-1}(0))))=(-2|c|^{2}-2\mbox{Re}(c\zeta))/(\overline{c\zeta}+1)$. Since $-2|c|^{2}-2\mbox{Re}(c\zeta)$ is a real number, $c\zeta$ is a real number. If $c \zeta > 0$, then $\tau(\varphi(\sigma(\tau^{-1}(0)))) < 0$. Therefore, $c \zeta < 0$ and $c \zeta=-|c|$.
Hence $\varphi(z)=(1-|c|)z/(cz+1)$.\hfill $\Box$ \\ \par

\section{Spectral radii of hyponormal weighted composition operators on $H^{2}$ and $A^{2}_{\alpha}$}
We know that if an operator $T$ is hyponormal, then $\|T\|=r(T)$. In this section, we investigate spectral radii of weighted composition operators. \\ \par

{\bf Lemma 3.1.} {\it Suppose that $\varphi$ is an analytic map of the unit disk into itself with Denjoy-Wolff point $\zeta \in \partial \mathbb{D}$. Then the essential spectral radius of $C_{\varphi}$ on $H^{2}$ or $A_{\alpha}^{2}$ is $\varphi'(\zeta)^{-\gamma/2}$.}\bigskip

{\bf Proof.} By using the general version of the Chain Rule given in \cite[Chapter 4, Exercise 10, p. 74]{sh}, we see that $\varphi_{n}'(\zeta)=\varphi'(\zeta)^{n}$. By \cite[Theorem 5.2]{km},
$r_{e,\gamma}(C_{\varphi})=\lim_{n \rightarrow \infty}\|C_{\varphi_{n}}\|_{e,\gamma}^{1/n}
 \geq \lim_{n \rightarrow \infty} \left(\frac{1}{|\varphi'_{n}(\zeta)|^{\gamma/2}} \right)^{1/n}
=\varphi'(\zeta)^{-\gamma/2}.$
  Also by \cite[Theorem 3.9]{cm1} and \cite[Theorem 6]{hu}, $\varphi'(\zeta)^{-\gamma/2}=r_{\gamma}(C_{\varphi}) \geq r_{e,\gamma}(C_{\varphi})$, so the result follows.\hfill $\Box$ \\ \par

Suppose that $\varphi$
is an analytic self-map of $\mathbb{D}$ and $\alpha$ is a complex
number of modulus 1. Since
Re$\left(\frac{\alpha+\varphi}{\alpha-\varphi}\right)$ is a
positive harmonic function on $\mathbb{D}$, there exists a finite
positive Borel measure $\mu_{\alpha}$ on $\partial \mathbb{D}$
such that
$\frac{1-|\varphi(z)|^{2}}{|\alpha-\varphi(z)|^{2}}=
\mbox{Re}\left(\frac{\alpha+\varphi(z)}{\alpha-\varphi(z)}\right)\\
=\int_{\partial \mathbb{D}}P_{z} d\mu_{\alpha}$ for each $z \in
\mathbb{D}$, where
$P_{z}(e^{i\theta})=(1-|z|^{2})/|e^{i\theta}-z|^{2}$ is the
Poisson kernel at $z$. The measures $\mu_{\alpha}$ are called the
Clark measures of $\varphi$. There is a unique pair  of measures
$\mu_\alpha^{ac}$ and $\mu_\alpha^s$ such that
$\mu_\alpha=\mu_\alpha^{ac}+\mu_\alpha^s$, where $\mu_\alpha^{ac}$
and $\mu_\alpha^s$ are the absolutely continuous and singular
parts with respect to Lebesgue measure, respectively. In particular, if $\varphi$ is a
linear-fractional non-automorphism such that
$\varphi(\zeta)=\eta$ for some $\zeta,\eta \in \partial
\mathbb{D}$, then $\mu_{\alpha}^{s}=0$ when $\alpha \neq \eta$ and
$\mu_{\eta}^{s}=|\varphi^{'}(\zeta)|^{-1}\delta_{\zeta}$, where $\delta_{\zeta}$ is the unit point mass at $\zeta$. We write
$E(\varphi)$ for the closure in $\partial \mathbb{D}$ of the
union of the closed supports of $\mu_{\alpha}^{s}$ as $\alpha$
ranges over the unit circle. We know that
$F(\varphi)\subseteq E(\varphi)$ (see \cite[p. 2919]{km}). For information about the Clark measures, see \cite{km}.\\ \par

{\bf Lemma 3.2.} {\it Let $\varphi$ be an analytic self-map of $\mathbb{D}$. Suppose that $\varphi \in A(\mathbb{D})$ and the set of points which $\varphi$ makes contact with $\partial \mathbb{D}$ is finite. Assume that there are a positive integer $n$ and $\zeta \in \partial\mathbb{D}$ such that $E(\varphi_{n})=\{\zeta\}$, where $\zeta$ is the Denjoy-Wolff point of $\varphi$. Let $\psi \in H^{\infty}$ be continuous at $\zeta$. Then
$$r_{\gamma}(C_{\psi,\varphi})=|\psi(\zeta)|\varphi'(\zeta)^{-\gamma/2}.$$
}\bigskip

{\bf Proof.} Since $\zeta$ is the Denjoy-Wolff point of $\varphi$, $F(\varphi_{n})=\{\zeta\}$. Also we have $\varphi_{n}'(\zeta)=\varphi'(\zeta)^{n}$.
Since $\sigma_{e,\gamma}(C_{\varphi_{n}})$ is a compact set, by Lemma 3.1, there is $\lambda \in \sigma_{e,\gamma}(C_{\varphi_{n}})$ such that $|\lambda|=r_{e,\gamma}(C_{\varphi_{n}})=\varphi'(\zeta)^{-n \gamma/2}$. Since $\varphi(\zeta)=\zeta$, \cite[Corollary 2.2]{kmm} and \cite[Proposition 2.3]{fash} imply that
$\sigma_{e,\gamma}(C^{n}_{\psi,\varphi})=\sigma_{e,\gamma}(T_{\psi \cdot \psi\circ\varphi...\psi\circ\varphi_{n-1}}C_{\varphi_{n}})
=\psi(\zeta)^{n}\sigma_{e,\gamma}(C_{\varphi_{n}}).$
We may now apply Lemma 2.1 and \cite[Lemma 5.1]{Bourdon1} to observe that  $|\mu|\leq |\psi(\zeta)|^{n} r_{\gamma}(C_{\varphi_{n}})$ for each $\mu$ in $\sigma_{p,\gamma}(C^{n}_{\psi,\varphi})$. Then by  \cite[Theorem 6]{hu} and \cite[Theorem 3.9]{cm1}, we have  for each $\mu \in \sigma_{p,\gamma}(C^{n}_{\psi,\varphi})$, $|\mu|\leq |\psi(\zeta)|^{n}\varphi'(\zeta)^{-n \gamma/2}$. By  \cite[Proposition 6.7, p. 210]{c1} and \cite[Proposition 4.4, p. 359]{c1},
$$\partial \sigma_{\gamma} (C^{n}_{\psi,\varphi})\subseteq \sigma_{ap,\gamma}(C^{n}_{\psi,\varphi})\subseteq \sigma_{p,\gamma}(C^{n}_{\psi,\varphi}) \cup \sigma_{e,\gamma}(C^{n}_{\psi,\varphi}).$$
Therefore, we can easily see that $r_{\gamma}(C^{n}_{\psi,\varphi})=|\psi(\zeta)|^{n}\varphi'(\zeta)^{-n \gamma/2}$. We have
$$r_{\gamma}(C_{\psi,\varphi})=\lim_{k\rightarrow\infty}\|C^{nk}_{\psi,\varphi}\|_{\gamma}^{1/(nk)}=\left(r_{\gamma}(C^{n}_{\psi,\varphi})\right)^{1/n}=|\psi(\zeta)|\varphi'(\zeta)
^{-\gamma/2}.$$
\hfill $\Box$ \\ \par

{\bf Corollary 3.3.} {\it Suppose that $\varphi$ is a parabolic non-automorphism with a fixed point $\zeta \in \partial\mathbb{D}$. Let $\psi \in H^{\infty}$ be continuous at $\zeta$. Then on $H^{2}$ and $A_{\alpha}^{2}$,
$$r_{\gamma}(C_{\psi,\varphi})=|\psi(\zeta)|.$$
}\bigskip

{\bf Proof.} Since $\varphi'(\zeta)=1$, the desired conclusion follows from Lemma 3.2.
\hfill $\Box$ \\ \par

In \cite[Theorem 21]{derek}, Cowen et al. showed that a hyponormal weighted composition operator on $H^{2}$ with a composition symbol $\varphi$ as a positive parabolic non-automorphism is automatically normal. In the following proposition, for  $\varphi$ which is a parabolic non-automorphism,  we find a necessary condition for $C_{\psi,\varphi}$ to be hyponormal on $H^{2}$ and $A^{2}_{\alpha}$. \\ \par

{\bf Proposition 3.4.} {\it Let $\varphi$ be a parabolic non-automorphism with a fixed point $\zeta \in \partial\mathbb{D}$. Suppose that $\psi \in H^{\infty}$ is continuous at $\zeta$. If $C_{\psi,\varphi}$ is hyponormal on $H^2$ or $A_{\alpha}^{2}$, then for each $w \in \mathbb{D}$,
\begin{equation}
|\psi(\zeta)|\geq|\psi(w)|\left(\frac{1-|w|^{2}}{1-|\varphi(w)|^{2}}\right)^{\gamma/2}.
 \end{equation}
}\bigskip

{\bf Proof.} Since $C_{\psi,\varphi}$ is hyponormal, by \cite[Proposition 4.6 p. 47]{c4} and Corollary 3.3, $\|C_{\psi,\varphi}\|_{\gamma}=|\psi(\zeta)|$. By Equation (1), we have
$$\|C_{\psi,\varphi}\|_{\gamma}=\|C_{\psi,\varphi}^{\ast}\|_{\gamma}
 \geq \|C_{\psi,\varphi}^{\ast}(K_{w}/\|K_{w}\|_{\gamma})\|_{\gamma}
=|\psi(w)|\left(\frac{1-|w|^{2}}{1-|\varphi(w)|^{2}}\right)^{\gamma/2},$$
as desired.
\hfill $\Box$ \\ \par

{\bf Example 3.5.}  Let $\varphi$ be a parabolic non-automorphism with fixed point 1. Assume that
$\psi_{1}(z)=\frac{2-z}{4}$ and $\psi_{2}(z)=-3z^2+2z+3.$
Setting $w=0$ and $\zeta=1$ in Equation (7), we get $C_{\psi_1,\varphi}$ and $C_{\psi_2,\varphi}$ are not hyponormal.
\bigskip





{\bf Proposition 3.6.} {\it Let $\varphi$ be an analytic self-map of $\mathbb{D}$
with $\varphi(0)=0$ and $\psi \in H^{\infty}$. Suppose there is a positive integer $n$  that $\{e^{i\theta}:|\varphi_{n}(e^{i\theta})|=1\}=\emptyset$. If $C_{\psi,\varphi}$ is hyponormal on $H^2$ or $A_{\alpha}^{2}$, then $\|C_{\psi,\varphi}\|_{\gamma}=r_{\gamma}(C_{\psi,\varphi})=|\psi(0)|$.
 }\bigskip

{\bf Proof.} Since $\overline{\varphi_{n}({\mathbb{D}})}\subseteq\mathbb{D}$, by \cite[p. 129]{cm1}, $C^{n}_{\psi,\varphi}=C_{\psi \cdot \psi\circ\varphi...\psi\circ\varphi_{n-1},\varphi_{n}}$ is compact. Therefore, $\sigma_{e,\gamma}(C^{n}_{\psi,\varphi})=\{0\}$ and so by the Spectral Mapping Theorem, $\sigma_{e,\gamma}(C_{\psi,\varphi})=\{0\}$. Schwarz's Lemma implies that 0 is the Denjoy-Wolff point. Invoking Lemma 2.1 and \cite[Lemma 5.1]{Bourdon1}, for each $\lambda \in \sigma_{p,\gamma}(C_{\psi,\varphi})$, $|\lambda|\leq|\psi(0)|r_{\gamma}(C_{\varphi})$. Hence by \cite[Theorem 6]{hu} and \cite[Theorem 3.9]{cm1}, for each $\lambda \in \sigma_{p,\gamma}(C_{\psi,\varphi})$, $|\lambda|\leq|\psi(0)|$. By \cite[Proposition 6.7, p. 210]{c1} and \cite[Proposition 4.4, p. 359]{c1},
$\partial\sigma_{\gamma}(C_{\psi,\varphi})\subseteq \sigma_{e,\gamma}(C_{\psi,\varphi})\cup \sigma_{p,\gamma}(C_{\psi,\varphi})$. Hence $r_{\gamma}(C_{\psi,\varphi})\leq|\psi(0)|$. We have
 $$\|C_{\psi,\varphi}\|_{\gamma}=\|C^{\ast}_{\psi,\varphi}\|_{\gamma}\geq\|C^{\ast}_{\psi,\varphi}K_{0}\|_{\gamma}=|\psi(0)|\|K_{\varphi(0)}\|_{\gamma}=|\psi(0)|.$$ Since $C_{\psi,\varphi}$ is hyponormal, $\|C_{\psi,\varphi}\|_{\gamma}=r_{\gamma}(C_{\psi,\varphi})=|\psi(0)|$.\hfill $\Box$ \\ \par

{\bf Proposition 3.7.} {\it Let $\varphi$ be analytic on $\mathbb{D}$ with $\varphi(\mathbb{D}) \subseteq \mathbb{D}$ and $\varphi(0)=0$. Assume that there is an integer $n$ such that $\{e^{i\theta}:|\varphi_{n}(e^{i\theta})|=1\}$ has only one element $\zeta$ which is a fixed point of $\varphi$ and $\zeta \in F(\varphi)$. Suppose that
$\varphi \in A(\mathbb{D})$ and $\psi \in H^{\infty}$ is continuous at $\zeta$. If $C_{\psi,\varphi}$ is hyponormal on $H^{2}$ or $A_{\alpha}^{2}$, then
$$\frac{|\psi(\zeta)|}{|\varphi'(\zeta)|^{\gamma/2}} \leq \|C_{\psi,\varphi}\|_{\gamma} \leq  \mbox{max} \{|\psi(\zeta)|,|\psi(0)|\}.$$
}\bigskip

{\bf Proof.} Assume that $\varphi(0)=0$ and $\varphi(\zeta)=\zeta$, where $\zeta \in \partial \mathbb{D}$.
By \cite[Corollary 2.2]{kmm} and \cite[Proposition 2.3]{fash},
$\sigma_{e,\gamma}(C_{\psi,\varphi}^{n})= \sigma_{e,\gamma}(C_{\psi \cdot \psi \circ \varphi ... \psi \circ  \varphi_{n-1},\varphi_{n}})\\=\psi(\zeta)^{n}\sigma_{e,\gamma}(C_{\varphi_{n}}).$
By \cite[Corollary 3.7]{cm1} and \cite[Lemma 2.3]{r}, we see that $r_{e,\gamma}(C_{\psi,\varphi}^{n}) \leq |\psi(\zeta)|^{n}r_{e,\gamma}(C_{\varphi_{n}}) \leq |\psi(\zeta)|^{n}\|C_{\varphi_{n}}\|_{\gamma} \leq |\psi(\zeta)|^{n}$. Hence by the Spectral Mapping Theorem, $r_{e,\gamma}(C_{\psi,\varphi}) \leq |\psi(\zeta)|$. As we saw in the proof of Proposition 3.6, \cite[Theorem 3.9]{cm1}, \cite[Theorem 6]{hu}, Lemma 2.1 and \cite[Lemma 5.1]{Bourdon1} imply that for each $\lambda \in \sigma_{p,\gamma}(C_{\psi,\varphi})$, $|\lambda|\leq|\psi(0)|$. Also by \cite[Proposition 6.7, p. 210]{c1} and \cite[Proposition 4.4, p. 359]{c1},
$\partial\sigma_{\gamma}(C_{\psi,\varphi})\subseteq \sigma_{e,\gamma}(C_{\psi,\varphi})\cup \sigma_{p,\gamma}(C_{\psi,\varphi}).$
Hence $r_{\gamma}(C_{\psi,\varphi}) \leq \mbox{max} \{|\psi(\zeta)|,|\psi(0)|\}$. By Equation (1), we see that for each $r <1$,
$$\|C_{\psi,\varphi}\|^{2}_{\gamma}
 = \|C_{\psi,\varphi}^{\ast}\|^{2}_{\gamma}
\geq  \left\|C_{\psi,\varphi}^{\ast}\frac{K_{r\zeta}}{\|K_{r \zeta}\|_{\gamma}}\right\|^{2}_{\gamma}
= |\psi(r\zeta)|^{2}\left(\frac{1-|r \zeta|^{2}}{1-|\varphi (r \zeta)|^{2}}\right)^{\gamma}.$$
Then by the Julia-Caratheodory Theorem (see \cite[Theorem 2.44, p. 51]{cm1}),
$$\|C_{\psi,\varphi}\|^{2}_{\gamma}
\geq  \lim_{r \rightarrow 1}|\psi(r\zeta)|^{2}\left(\frac{1+|r \zeta|}{1+|\varphi(r \zeta)|}\right)^{\gamma}\left(\frac{1-|r \zeta|}{1-|\varphi(r\zeta)|}\right)^{\gamma}
=|\psi(\zeta)|^2\frac{1}{|\varphi'(\zeta)|^{\gamma}}.$$
Since $\|C_{\psi,\varphi}\|_{\gamma}=r_{\gamma}(C_{\psi,\varphi})$, the result follows.\hfill $\Box$ \\ \par
Let $H$ be a Hilbert space of functions analytic on the unit disk. If the monomials $1, z, z^{2},...$ are an orthogonal set of non-zero vectors with dense span in $H$, then $H$ is called a weighted Hardy space. We will assume that the norm satisfies the normalization $\|1\|=1$. The weight sequence for a weighted Hardy space $H$ is defined to be $\beta(n)=\|z^{n}\|$. The weighted Hardy space with weight sequence $\beta(n)$ will be denoted  $H^{2}(\beta)$. The norm on $H^{2}(\beta)$ is given by
$\left\|\sum_{j=0}^{\infty}a_{j}z^{j}\right\|^{2}=\sum_{j=0}^{\infty}|a_{j}|^{2}\beta(j)^{2}.$
We know that $H^{2}$ and $A^{2}_{\alpha}$ are weighted Hardy spaces (see \cite{cm1}).\par
Suppose that $\varphi$, not the identity and not an elliptic automorphism of $\mathbb{D}$, is an analytic map of the unit disk into itself with $\varphi(0)=0$.
In the following proposition, we show that for a hyponormal weighted composition operator $C_{\psi,\varphi}$, if $r_{e,\gamma}(C_{\psi,\varphi}) \leq |\psi(0)|$, then $\psi$ is constant and $C_{\varphi}$ is hyponormal. \\ \par

{\bf Proposition 3.8.} {\it Suppose that $\varphi$, not the identity and not an elliptic automorphism of $\mathbb{D}$, is an analytic map of the unit disk into itself with $\varphi(0)=0$. Assume that $\psi \in H^{\infty}$ and for each $\lambda \in \sigma_{e,\gamma}(C_{\psi,\varphi})$, $|\lambda|\leq |\psi(0)|$. If $C_{\psi,\varphi}$ is hyponormal on $H^2$ or $A_{\alpha}^{2}$,
then $\psi$ is constant, $C_{\varphi}$ is hyponormal and $\|C_{\psi,\varphi}\|_{\gamma}=r_{\gamma}(C_{\psi,\varphi})=|\psi(0)|$.}\bigskip

{\bf Proof.} By \cite[Lemma 5.1]{Bourdon1}, Lemma 2.1, \cite[Theorem 3.9]{cm1} and \cite[Theorem 6]{hu}, for each $\lambda \in \sigma_{p,\gamma}(C_{\psi,\varphi})$, $|\lambda|\leq |\psi(0)|$. \cite[Proposition 6.7, p. 210]{c1} and \cite[Proposition 4.4, p. 359]{c1} imply that $r_{\gamma}(C_{\psi,\varphi})\leq |\psi(0)|$. Since $C_{\psi,\varphi}$ is hyponormal,
$|\psi(0)|\geq\|C_{\psi,\varphi}\|_{\gamma}\geq\|C_{\psi,\varphi}1\|_{\gamma}=\|\psi\|_{\gamma}.$
Since $H^{2}$ and $A_{\alpha}^{2}$ are weighted Hardy spaces, it is easy to see that $\|\psi\|_{\gamma} \geq |\psi(0)|$.
Hence $\psi$ is constant and $\|C_{\psi,\varphi}\|_{\gamma}=r_{\gamma}(C_{\psi,\varphi})=|\psi(0)|$.\hfill $\Box$ \\ \par


By the definition of hyponormality, we can see that if a hyponormal operator is unitarily equivalent to another operator, then that other operator is also hyponormal; we use this fact in the proof of the following two theorems. \par




Suppose that $\varphi$, not the identity and not an elliptic automorphism of $\mathbb{D}$, is an analytic map of the unit disk into itself. In the following theorem, we see that if $C_{\psi,\varphi}$ is hyponormal, when $\varphi(p)=p$ for some $p \in \mathbb{D}$ and $r_{e,\gamma}(C_{\psi,\varphi}) \leq |\psi(p)|$, then the function $\psi$ has a simple linear-fractional form that is the same as what Bourdon et al. and Le found in  \cite[Theorem 10]{Bourdon} and \cite[Theorem 4.3]{l}. Also Theorem 3.9 is an extension of Proposition 3.8.\\ \par

{\bf Theorem 3.9.} {\it Suppose that $\varphi$, not the identity and not an elliptic automorphism of $\mathbb{D}$, is an analytic map of the unit disk into itself with $\varphi(p)=p$, where $p \in \mathbb{D}$.  Assume that $\psi \in H^{\infty}$. Suppose that  for each $\lambda \in \sigma_{e,\gamma}(C_{\psi,\varphi})$, $|\lambda|\leq |\psi(p)|$. If $C_{\psi,\varphi}$ is hyponormal on $H^2$ or $A_{\alpha}^{2}$,
then $\|C_{\psi,\varphi}\|_{\gamma}=r_{\gamma}(C_{\psi,\varphi})=|\psi(p)|$ and
 $\psi=\psi(p)\frac{K_{p}}{K_{p} \circ \varphi}$.}\bigskip

{\bf Proof.} Let $C_{\psi,\varphi}$ be hyponormal. Suppose that $\psi_{p}=K_{p}/\|K_{p}\|_{\gamma}$ and $\alpha_{p}(z)=(p-z)/(1-\overline{p}z)$. By \cite[Theorem 6]{Bourdon} and \cite[Corollary 3.6]{l}, we have
\begin{eqnarray}
H:=C^{\ast}_{\psi_{p},\alpha_p}C_{\psi,\varphi}C_{\psi_{p},\alpha_p}
\end{eqnarray}
is hyponormal. Let $\sigma$ be the Krein adjoint of $\alpha_{p}$ and $g$ and $h$ be the Cowen auxiliary functions for $\alpha_{p}$. Since $T^{\ast}_{h}T^{\ast}_{\psi_{p}}=(T_{\psi_{p}h})^{\ast}=T_{1/\|K_{p}\|_{\gamma}}$ and $C_{\alpha_{p}}^{\ast}=T_{g}C_{\alpha_{p}}T_{h}^{\ast}$, by \cite[Remark 2.1(a)]{fash}, we see that
\begin{eqnarray}
  H&=&T_{g}C_{\alpha_{p}}T^{\ast}_{h}T^{\ast}_{\psi_{p}}T_{\psi}C_{\varphi}T_{\psi_{p}}C_{\alpha_p}\nonumber\\
  &=&\frac{1}{\|K_{p}\|_{\gamma}}T_{g}T_{\psi \circ \alpha_{p}}T_{\psi_{p}\circ \varphi\circ\alpha_{p}}C_{\alpha_{p}}C_{\varphi}C_{\alpha_{p}}\nonumber\\
  &=&C_{q,\alpha_{p}\circ \varphi\circ\alpha_{p}},
\end{eqnarray}
where $q=(g \cdot \psi \circ \alpha_{p} \cdot \psi_{p} \circ \varphi \circ \alpha_{p})/\|K_{p}\|_{\gamma}$ (see the proof of \cite[Theorem 10]{Bourdon}).
 Since $H$ is unitarily equivalent to $C_{\psi,\varphi}$, $\sigma_{e,\gamma}(C_{\psi,\varphi})=\sigma_{e,\gamma}(H)$. Since $q(0)=\psi(p)$, for each $\lambda \in \sigma_{e,\gamma}(C_{q,\alpha_{p} \circ \varphi \circ \alpha_{p}})$,  $|\lambda|\leq |q(0)|$. We may apply Proposition 3.8 to conclude that $\|C_{q,\alpha_{p} \circ \varphi \circ \alpha_{p}}\|_{\gamma}=r_{\gamma}(C_{q,\alpha_{p} \circ \varphi \circ \alpha_{p}})=|q(0)|=|\psi(p)|$. Since  $r_{\gamma}(C_{q,\alpha_{p} \circ \varphi \circ \alpha_{p}})=r_{\gamma}(C_{\psi,\varphi})$ and $C_{\psi,\varphi}$ is hyponormal, $\|C_{\psi,\varphi}\|_{\gamma}=r_{\gamma}(C_{\psi,\varphi})=|\psi(p)|$, as desired. Also Proposition 3.8 implies that $q$ is constant. Since $q \equiv \psi(p)$, $g \cdot \psi\circ\alpha_{p} \cdot \psi_{p}\circ\varphi\circ\alpha_{p} \equiv \|K_{p}\|_{\gamma}\psi(p)$ on $\mathbb{D}$. Then
$g\circ\alpha_{p} \cdot \psi \cdot \psi_{p}\circ\varphi \equiv \|K_{p}\|_{\gamma}\psi(p)$ on $\mathbb{D}$. It follows that
$\psi=\|K_{p}\|_{\gamma}\psi(p)/(g\circ\alpha_{p} \cdot \psi_{p}\circ\varphi).$
Observe
$g\circ\alpha_{p}=\|K_{p}\|_{\gamma}^{2}/K_{p}$
and
 $\psi_{p}\circ\varphi=K_{p}\circ\varphi/\|K_{p}\|_{\gamma}$.
Hence
$\psi=\psi(p)\frac{K_{p}}{K_{p}\circ\varphi}$.
\hfill $\Box$ \\ \par

{\bf Corollary 3.10.} {\it Suppose that $\varphi$ is analytic in a neighborhood of the closed unit disk. Let there exist $p \in \mathbb{D}$ such that $\varphi(p)=p$. Assume that there is a positive integer $n$ such that $\{e^{i\theta}:|\varphi_{n}(e^{i\theta})|=1\}$ has only one element $\zeta$ which is a fixed point of $\varphi$ and $\zeta \in F(\varphi)$.  Suppose that $\psi \in H^{\infty}$  is continuous at $\zeta$ and $|\psi(\zeta)|  \leq |\psi(p)|$. If $C_{\psi,\varphi}$ is hyponormal on $H^{2}$, then $\psi = \psi(p)\frac{K_{p}}{K_{p} \circ \varphi}$ and $\|C_{\psi,\varphi}\|_{1}=|\psi(p)|$ .}\bigskip

{\bf Proof.} By \cite[Corollary 2.2]{kmm},
$\sigma_{e,1}(C_{\psi,\varphi}^{n})=\sigma_{e,1}(T_{\psi \cdot \psi \circ \varphi...\psi \circ \varphi_{n-1}C_{\varphi_{n}}})=\psi(\zeta)^{n}\sigma_{e,1}(C_{\varphi_{n}}).$
Also \cite[Exercise 3.2.5]{cm1} and the general version of the Chain Rule given in \cite[Chapter 4, Exercise 10, p. 74]{sh} imply that $|\lambda| \leq |\psi(\zeta)|^{n}r_{e,1}(C_{\varphi_{n}}) \leq |\psi(\zeta)|^{n}|\varphi'(\zeta)|^{-n/2} \leq |\psi(\zeta)|^{n}$ for each $\lambda \in \sigma_{e,1}(C_{\psi,\varphi}^{n})$.  By  the Spectral Mapping Theorem, we can see that for each $t \in \sigma_{e,1}(C_{\psi,\varphi})$, $|t| \leq |\psi(\zeta)|$. The result follows from Theorem 3.9.\hfill $\Box$ \\ \par

One example of $\varphi$ for Corollary 3.10 is a hyperbolic non-automorphism with a fixed point in $\mathbb{D}$ and another fixed point in $\partial \mathbb{D}$.\\ \par

{\bf Theorem 3.11.} {\it Let  $\varphi \in A(\mathbb{D})$ with $\varphi(\mathbb{D})\subseteq \mathbb{D}$. Assume $p$ and $\zeta$ satisfy the hypotheses of the second and third sentences of Corollary 3.10.
 Suppose $\psi \in H^{\infty}$ is continuous at $\zeta$. If $C_{\psi,\varphi}$ is hyponormal on $H^{2}$ or $A_{\alpha}^{2}$, then
$$ \frac{\mu}{|\varphi'(\zeta)|^{\gamma/2}}\leq \|C_{\psi,\varphi}\|_{\gamma}\leq \max\{\mu,|\psi(p)|\},$$
where $\mu=|\psi(\zeta)K_{p}(\alpha_{p}(\zeta))K_{p}(\zeta)|/ \|K_p\|^2_{\gamma}$ with $\alpha_{p}(z)=(p-z)/(1-\bar{p}z)$.}\bigskip

{\bf Proof.} Suppose that $C_{\psi,\varphi}$ is hyponormal. Let $H$ be as in Equation (8). We see that the map $\alpha_{p} \circ \varphi \circ \alpha_{p}$ fixes $\alpha_p(\zeta)$. Since $\alpha_{p}^{-1}=\alpha_{p}$, by \cite[Corollary 7.6, p. 99]{c3}, $(\alpha_{p} \circ \varphi \circ \alpha_{p})'(\alpha_{p}(\zeta))=\varphi'(\zeta)$. By Proposition 3.7, we have
\begin{eqnarray}
\frac{|q(\alpha_p(\zeta))|}{|\varphi'(\zeta)|^{\gamma/2}}\leq \|C_{q,\alpha_{p} \circ \varphi \circ \alpha_{p}}\|_{\gamma}\leq \max\{|q(\alpha_{p}(\zeta))|,|q(0)|\},
\end{eqnarray}
where $q=(g \cdot \psi \circ \alpha_p \cdot \psi_p \circ \varphi \circ \alpha_{p})/\|K_{p}\|_{\gamma}$ with $g$  the Cowen auxiliary function for $\alpha_{p}$ and $\psi_{p}=K_{p}/\|K_{p}\|_{\gamma}$. We have $q(0)=\psi(p)$. Moreover,
$$q(\alpha_p(\zeta))=\frac{\psi(\zeta)K_{p}(\alpha_p(\zeta))K_{p}(\zeta)}{\|K_{p}\|_{\gamma}^2}.$$
  Since $H$ is unitarily equivalent to $C_{\psi,\varphi}$, $r_{\gamma}(H)=r_{\gamma}(C_{\psi,\varphi})$. Since $C_{\psi,\varphi}$ and $H$ are hyponormal, by Equations (9) and (10), the result follows.\hfill $\Box$ \\ \par

\textbf{Acknowledgments}
\bigskip

The authors really appreciate professor  Barbara D. MacCluer  of the University of Virginia for her attention and practical help. Also the authors would like to thank professor Derek Thompson of Taylor University for pointing out some useful remarks.

\bigskip
{
M. Fatehi, Department of Mathematics, Shiraz Branch, Islamic Azad
University, Shiraz, Iran. \par E-mail: fatehimahsa@yahoo.com \par
M. Haji Shaabani, Department of Mathematics, Shiraz University of Technology, P. O. Box 71555-313, Shiraz, Iran.\par E-mail: shaabani@sutech.ac.ir\par}

\end{document}